\documentclass[11pt,amssym,twoside]{article}
\usepackage{amssymb}
\newtheorem{thm}{Theorem}[section]
\newtheorem{cor}[thm]{Corollary}
\newtheorem{lem}[thm]{Lemma}
\newtheorem{rmk}[thm]{Remark}
\newtheorem{prop}[thm]{Proposition}
\newtheorem{problem}[thm]{Problem}
\newtheorem{conj}[thm]{Conjecture}

\newtheorem{defn}[thm]{Definition}

\newcommand{\Real}{\mathbb R}
\newcommand{\Cplx}{\mathbb C}

\newcommand{\Z}{\mathbb Z}
\newcommand{\F}{\mathbb F}
\newcommand{\G}{\mathbb G}
\newcommand{\N}{\mathbb N}
\newcommand{\Proj}{\mathbb P}
\newcommand{\Q}{\mathbb Q}
\begin{document}

\title{Self-Similar Fractals and Arithmetic Dynamics}%
\author{Arash Rastegar}%
%\address{Sharif University of Technology}%
%\email{rastegar@sharif.edu}%

%\thanks{}%
%\subjclass{}%
%\keywords{}%

%\date{}%
%\dedicatory{}
%\commby{}%
% ----------------------------------------------------------------
\maketitle
\begin{abstract}
The concept of self-similarity on subsets of algebraic varieties
is defined by considering algebraic endomorphisms of the variety
as `similarity' maps. Self-similar fractals are subsets of algebraic varieties
which can be written as a finite and disjoint union of
`similar' copies. Fractals provide a framework in which, one can
unite some results and conjectures in Diophantine geometry. We
define a well-behaved notion of dimension for self-similar fractals. We also
prove a fractal version of Roth's theorem for algebraic points on
a variety approximated by elements of a fractal subset. As a
consequence, we get a fractal version of Siegel's theorem on finiteness of integral points
on hyperbolic curves and a fractal version of Falting's theorem on Diophantine approximation on abelian varieties.
\end{abstract}
\section*{Introduction}

Self-similar fractals are very basic geometric objects which
presumably could have been defined as early as Euclid. By
self-similar fractals, we mean objects which are (almost) disjoint union of pieces
`similar' to the whole object. In Euclidean context, one can
think of Euclidean plane as the ambient space and Euclidean
similarities as `similarity' maps. There are several interesting
examples of such fractals in the literature. Sierpinski carpet,
Koch snowflake, and Cantor set are among the typical examples of
Euclidean fractals. In a more modern geometric context, the
ambient space of an affine fractal could be a real vector space,
and `similarity' maps could be chosen to be affine maps, which are
usually assumed to be distance decreasing.

In the algebraic context, ambient space of an affine fractal
could be a vector space over arbitrary field and polynomial
self-maps of the vector space with coefficients in the base field
could be taken as `similarity' maps. Ideals in the ring of
integers of a number field are examples of affine fractals in
this context. Self-similar fractals in a ring could be much more
complicated. For example, integers missing a number of digits in
their decimal expansion form a fractal. The algebraic concept of
self-similar fractals could also be extended to subsets of
algebraic varieties, if we take algebraic endomorphisms as
`similarity' maps.

In this paper, we assume that a fractal is a finite union of its
similar images except for finitely many points 
and a self-similar fractal has the extra condition that
this fractal is a finite union of its
similar images and its similar images are disjoint or at most with finite intersection. For
example, rational points on a projective space could be thought
of as a self-similar set, but not as a fractal, since it is union 
of infinitely many similar copies of itself.

The first important question about self-similar fractals is how to define
their dimension. One can introduce a notion of dimension which is
independent of the representation of the self-similar fractal as union of
similar images. We use arithmetic height-functions to introduce
such a concept of dimension for self-similar fractals. In fact, this notion of
fractal-dimension turns out to be related to the growth of the
number of points of bounded height in our fractal. This way, we
recover some classical computations in this direction.

One can think of Diophantine approximation of algebraic points by
a fractal whose elements are algebraic over $\Q$. Self-similarity
of fractals imply a strong version of Roth's theorem in this case.

One shall note that, fractals are not necessarily dense in the
ambient space with respect to complex topology. Therefore, such
approximation theorems are only interesting if we are
approximating a limiting point with respect to some Riemannian
metric. 

As a reward, we get fractal versions of Siegel's theorem on finiteness of integral points
and Falting's theorem on diophantine approximation on abelian varieties. 
Here are special cases, which could be formulated without any reference to fractals.
We have treated these special cases separatedly in [Ras1] and [Ras2]:

\begin{thm}
Let $X$ be an affine open subcurve of a connected smooth projectuve curve of genus $\geq 1$ defined
over $\mathbb{C}$  in the ambient affine space $\mathbb{A}^n(\mathbb{C})$ and let $F\subset \mathbb{A}^n(\mathbb{C})$
denote any finitely generated subgroup of $\mathbb{C}^n$ . Then $X(K )\cap F$ is finite.
\end{thm}

This implies that Siegel's theorem is an algebro-geometric fact, not an Arithmetic one.

\begin{thm}
Let $A$ be an abelian variety defined over a finitely generated subfield $K$ of $\mathbb{C}$. Let
$E$ is a geometrically irreducible subvariety of $A$ defined over $K$ and $F$ be a finitely generated subgroup of 
$A(K)$.                                                                                                                                                                                                                                                                       
Let $w$ be a valuation on $K$ and $H(x)$ a height function on $K$ coming from a choice of projective model for $K$
over the algebraic closure of $\mathbb{Q}$ in $K$.
If $d_w(x,E)$ denotes the $w$-adic distance from $x$ to $E$, and $\kappa$ and $c$ are positive constants,
then, there are only finitely many points in $F$ satisfying the following inequality
$$
d_w(x,E)< cH(x)^{-\kappa}.
$$
\end{thm}

This, in turn, implies that Faltings' theorem on Diophantine approximation on abelian varieties is also 
an algebro-geometric fact, not an Arithmetic one.

There are quite a few classical objects in arithmetic geometry
which can be considered as self-similar fractals. For example, for an abelian
variety $A$ defined over a number-field as ambient space, the set
of rational points $A(\Q)$ or any finitely generated subgroup of
$A(\bar {\Q})$ and the set of torsion points $A^{tor}$ can be
thought of as self-similar fractals with respect to endomorphisms of $A$.

In fact, fractals provide a common framework in which similar
theorems about objects in arithmetic geometry could be united in a
single context. For example, similarity between Manin-Mumford
conjecture on torsion points on an abelian variety which was
proved by Raynaud [Ray], and Lang's conjecture on finitely
generated subgroups of rational points on an abelian variety
which was proved by Faltings [Fal], made us propose the following general
conjecture about fractals:
\begin{conj}
Let $V$ be an irreducible variety defined over a finitely
generated field $K$ and let $F\subset V(K)$ denote a fractal
in $V$. Then, for any reduced subscheme $Z$ of $V$ defined over
$K$ the Zariski closure of $Z(K)\cap F$ is union of finitely
many points and finitely many components $B_j$ such that $B_j(K)\cap F$ 
is a fractal in $B_j$ with respect to some of the same self-similarity maps for each $j$.
\end{conj}
A generalized version of Lang's conjecture is covered by the above
conjecture. Some of our results in this paper also can be
considered as its special cases. Detailed evidences are presented
in the final section. We will also present a conjecture extending the above covering Andre-Oort conjecture,
proved by Pila and Tsimerman.

% ---------------------------------------------------------------
\section{Fractals in $\Z$}

The idea of considering fractal subsets of $\Z$ is due to O.
Naghshineh who proposed the following problem for "International
Mathematics Olympiad" held in Scotland in July 2002.
\begin{problem}
Let $F$ be an infinite subset of $\Z$ such that $F=\bigcup_{i=1}^n
a_i.F+b_i$ for integers $a_i$ and $b_i$ where $a_i.F+b_i$ and
$a_j.F+b_j$ are disjoint for $i\neq j$ and $|a_i|>1$ for each $i$.
Prove that
$$
\sum_{i=1}^n {1\over |a_i|}\leq 1.
$$
\end{problem}

In [Na], he explains his ideas about fractals in $\Z$ and suggests
how to define their dimension and how to prove this notion is
independent of the choice of self-similarity maps. His
suggestions are carried out by H. Mahdavifar. In this section, we
present their results and ideas.

\begin{defn}
Let $\phi_i:\Z \to \Z$ for $ i=1$ to $n$ denote linear maps of
the form $\phi_i(x)=a_i.x+b_i$ where $a_i$ and $b_i$ are integers
with $|a_i|>1$. A subset $F\subseteq \Z$ is called a self-similar fractal with
respect to $\phi_i$ if $F$ is disjoint union of its images under
the linear map $\phi_i$. In this case, we write $F=\sqcup_i
\phi_i(F)$ and define dimension of $F$ to be the real number $s$
such that
$$
\sum_{i=1}^n |a_i|^{-s}=1.
$$
\end{defn}

The basic example for self-similar fractals in $\Z$ is the set of integers
which miss a number of digits in their decimal expansion. This
definition of dimension is motivated by the notion of box
dimension for fractals on real vector spaces, which coincides
with Hausdorff dimension [Falc]. The challenge is to prove that,
this notion of dimension is independent of all the choices made,
and depends only on the self-similar fractal itself as a subset of $\Z$. Also,
smaller self-similar fractals should have smaller dimension. Having this
proven, it is easy to solve the above IMO problem. Note that $\Z$
is a self-similar fractal of dimension one. A self-similar fractal $F\subseteq \Z$ is of
dimension $\leq 1$ which solves the problem.

\begin{thm}
Let $F\subseteq \Z$ satisfy $F\subseteq\cup_i \phi_i(F)$ where
$\phi_i$ are as above. If $s$ is a real number such that $\sum_i
|a_i|^{-s}<1$ then the number of elements of $F$ in the ball
$B(x)$ is bounded above by $cx^s$ for some constant $c$ and for
large $x$.
\end{thm}
\textbf{Proof:} Let $F_i=\phi_i(F)$, and let $N(x)$ and $N_i(x)$
denote the number of elements of $F$ and $F_i$ in the ball
$B(x)$, respectively. We have
$$
N(x)\leq \sum_i N_i(x)
$$
and since for $f \in F_i$ and $\phi_i^{-1}(f)\in F$ we have
$|\phi_i^{-1}(f)|\leq (|f|+|b_i|)/|a_i|$ \, we can write
$$
N_i(x) \leq N({{x+|b_i|} \over|a_i|})
$$
If we let $t=Max_i\{|b_i|/|a_i|\}$ then we get the following
estimate
$$
N(x)\leq \sum_i N({x \over |a_i|}+t)
$$
We define a function $h:[1,\infty ] \to \Real$ by
$h(x)=x^{-s}N(x)$ and we shall show that $h$ is a bounded
function. The above estimate will have the form
$$
h(x)\leq \sum_i ({1\over |a_i|}+{t \over x})^s h({x\over |a_i|}+t)
$$
There exists a constant $M$ such that for $x>M$ we have $(x/
|a_i|)+t<x$ for all $i$ and
$$
\sum_i ({1\over |a_i|}+{t \over x})^s<1
$$
Now, assume $|a_1|\leq ...\leq |a_n|$ and define $x_0=|a_n|(M-t)$
and $x_j=|a_1|(x_{j-1}-t)$ for $j\geq 1$. Then $x_j$ is an
unbounded decreasing sequence. The function $h$ is bounded on
$[M,x_0]$ and we inductively show that it has the same bound on
$[x_j,x_{j+1}]$: for if $x\in [x_j,x_{j+1}]$ then $(x/
|a_i|)+t\in [(x_j/ |a_i|)+t,x-{j+1}/ |a_i|)+t] \subset [M, x_j]$
and if by induction hypothesis we have $h(x/ |a_i|)+t)\leq c$ for
all $i$ then
$$
h(x)\leq \sum_i ({1\over |a_i|}+{t \over x})^s h({x\over
|a_i|}+t) < c\sum_i ({1\over |a_i|}+{t \over x})^s < c
$$
It remains to notice that $h$ is also bounded on $[1,M]$.
$\square$

\begin{thm}
Let $F\subseteq \Z$ satisfy $F\supseteq\sqcup_i \phi_i(F)$ where
$\phi_i$ are as above. If $r$ is a real number such that $\sum_i
Norm(a_i)^{-r}>1$ then the number of elements of $F$ in the ball
$B(x)$ is bounded below by $cx^r$ for some constant $c$ and for
large $x$.
\end{thm}
\textbf{Proof:} We use the notation in the proof of the previous
lemma. Since for $f \in F_i$ and $\phi_i^{-1}(f)\in F$ we have
$|\phi_i^{-1}(f)|\geq (|f|-|b_i|)/|a_i|$ \, and we get
$$
N_i(x) \geq N({{x-|b_i|} \over|a_i|})\geq N({x \over |a_i|}-t)
$$
where $t=Max_i\{|b_i|)/|a_i|\}$. Now, it remains to show that
$h:[1,\infty ] \to \Real$ defined by $h(x)=x^{-r}N(x)$ is bounded
below, which can be proved along the same line as the previous
lemma. $\square$

\begin{prop} Let $F_1\subseteq F_2 \subseteq \Z$ be fractals. Then the
notion of fractal dimension is well-defined and $dim(F_1) \leq
dim(F_2)$.
\end{prop}
\textbf{Proof:} Suppose $F=\sqcup_i \phi_i(F)=\sqcup_j \psi_j(F)$
where $\phi_i$ and $\psi_i$ are linear functions
$\phi_i(x)=a_i.x+b_i$ and $\psi_j(x)=c_j.x+d_j$. Assume $\sum_i
|a_i|^{-\alpha}=1$ and $\sum_i |c_j|^{-\beta}=1$. We must show
that $\alpha=\beta$. Suppose $\alpha <\beta$. Insert real numbers
$\alpha <s<r<\beta$. Since $F\subset\cup_i \phi_i(F)$ and $\sum_i
|a_i|^{-s}<1$, we get $N(x)\leq cx^s$ for large $x$ and since
$F\supseteq\sqcup_i \psi_j(F)$ and $\sum_i |c_j|^{-r}>1$, we get
$N(x)\geq cx^r$ for large $x$ which is a contradiction. Thus
$\alpha=\beta$.

Now, for fractals $F_1\subseteq F_2$ suppose that $F_1=\sqcup_i
\phi_i(F)$ and $F_2=\sqcup_j \psi_j(F)$ where $\phi_i$ and
$\psi_i$ functions as above, and let $\sum_i |a_i|^{-\alpha}=1$
and $\sum_i |c_j|^{-\beta}=1$. We must show that
$\alpha\leq\beta$. Suppose $\alpha>\beta$ and insert real numbers
$\alpha>r>s>\beta$. Then one can get a contradiction as above.
$\square$

Naghshineh and Mahdavifar also suggest that the same calculations
work for $\Z[i]$ if we use norm of a complex number instead of
absolute value for a real number. The same arguments indicates
that, the notion of dimension of a fractal is linked to asymptotic
behavior of the number of points of bounded norm.

% ---------------------------------------------------------------
\section{Affine fractals}

It would be more convenient for the reader, if we formulate the
most general form of an affine fractal, and then treat special
cases.

\begin{defn}
Let $X$ be an affine algebraic variety defined over a finitely
generated field $K$, and let $f_i$ for $i=1$ to $n$, denote
polynomial endomorphisms of $X$ of degrees $\geq 1$ with
coefficients in $K$. A subset $F\subset X(K )$ is called an
affine self-similar fractal with respect to $f_1,...,f_n$ if $F$ is almost
disjoint union of its images under the polynomial endomorphisms
$f_i$ for $i=1$ to $n$, in which case, by abuse of notation, we
write $F=\sqcup_i f_i(F)$. An affine fractal in $X$ is a subset
which is affine fractal with respect to some polynomial
endomorphisms $f_1,...,f_n$. Note that such a representation is
not unique. In case we only have $F=\cup_i f_i(F)$ except for finitely many points of $F$ which are outside 
$\cup_i f_i(F)$ we simply call $F$ a fractal. In case $f_i$ are height increasing,
this will always be equal to forward orbit with respect to $f_i$ 
of finitely many points.
\end{defn}
$\textbf {1.}$ Let $K$ be a number field and let $O_K$ denote its
ring of integers. One can take $O_K$ as ambient space and
polynomial maps $\phi_i:O_K \to O_K$ with coefficients in $O_K$ as
self-similarities. Let $a_i$ denote the leading coefficient of
$\phi_i$, and $n_i$ denote the degree of $\phi_i$. Fix an
embedding $\rho :K\hookrightarrow \mathbb C$. Assume
$Norm_{\rho}(a_i)>1$ in case $\phi_i$ is linear. Let $F\subseteq
O_K$ be an affine fractal with respect to $\phi_i$ for $i=1$ to
$n$. One can define the fractal-dimension of $F$ to be the real
number $s$ for which
$$
\sum_{i=1}^n Norm(a_i)^{-{s \over n_i}}=1.
$$
Arguments of the previous section hold almost line by line, if one
replaces the absolute value of an integer with the product of
various archemidian norms of an algebraic integer in $O_K$.
Therefore, we have the following result:

\begin{prop} 
the above notion of dimension for affine fractals in
$O_K$ is well-defined and well-behaved with respect to inclusion
of affine fractals, i.e. dimension of an affine fractal is
independent of the choice of self-similarities and compatible
with inclusion of fractal subsets. 
\end{prop}

for a fractal generated by finitely many points we have proved the following:
\begin{prop}
Let $F$ be a fractal with respect to $f_i$ as above, the number of points of norm bounded by $X$
is $O(X^s)$ where $s$ is determined by 
$$
\sum_{i=1}^n Norm(a_i)^{-{s \over n_i}}=1.
$$
\end{prop}

$\textbf {2.}$ Start from a
linear semi-simple algebraic group $G$ and a rational
representation $\rho :G\to GL(W_{\Q})$ defined over $\Q$. Let
$w_0\in W_{\Q}$ be a point whose orbit $V=w_0\rho (G)$ is Zariski
closed. Then the stabilizer $H\subset G$ of $w_0$ is reductive
and $V$ is isomorphic to $H\setminus G$. By a theorem of
Borel-Harish-Chandra $V(\Z)$ breaks up to finitely many $G(\Z)$
orbits [Bo-HC]. Thus the points of $V(\Z)$ are parametrized by
cosets of $G(\Z)$. Fix an orbit $w_0G(\Z)$ with $w_0$ in $G(\Z)$.
Then the stabilizer of $w_0$ is $H(\Z)=H\cap G(\Z)$.

The additive structure of $G$ allows one to define self-similar
subsets of $V(\Z)$ and study their asymptotic behavior using the
idea of fractal dimension. For example, one can define
self-similarities to be maps $\phi:V(\Z)\to V(\Z)$ of the form
$$
\phi(\omega_0\gamma)=\omega_0([n]\gamma+g_0)
$$
where $[n]$ denotes multiplication by $n$ in $G(\Z)$ and $g_0$ is
an element in $G(\Z)$. These similarity maps are expansive if
$n>1$ and lead to a notion of dimension for self-similar fractals in $V(\Z)$.
upper bound similar to above holds for fractals in $V(\Z)$.

Duke-Rudnick-Sarnak [D-R-S] putting some extra technical
assumptions, have determined the asymptotic behavior of
$$
N(V(\Z),x)=\sharp\{\gamma \in H(\Z)\setminus G(\Z):
||w_0\gamma||\leq x\}.
$$
They prove that there are constants $a\geq0 ,b>0$ and $c>0$ such
that
$$
N(V(\Z),x)\sim cx^a(log x)^b.
$$
Note that, the whole set $V(\Z)$ could not be a fractal, since the
asymptotic behavior of its points is not polynomial. \\ $\textbf
{3.}$ Here is an example of an affine self-similar fractal with respect to
nonlinear polynomial maps. The subset
$$
\{(2^i,2^j)\in \Q^2 |i,j\in \Z \}
$$
is an affine self-similar fractal with respect to $f_1(x_1,x_2)=(x_1^2,x_2^2)$,
$f_2(x_1,x_2)=(2x_1^2,x_2^2)$, $f_3(x_1,x_2)=(x_1^2,2x_2^2)$ and
$f_4(x_1,x_2)=(2x_1^2,2x_2^2)$. Notice that, after
projectivization, we still get a self-similar set in the
projective line $\Proj^1(\Q)$. The subset
$$
\{(2^i;2^j)\in \Proj^1(\Q)|i,j\in \N \cup \{ 0\}\}
$$
is a self-similar fractal with respect to $f_1(x_1;x_2)=(x_1^2;x_2^2)$ and
$f_2(x_1;x_2)=(2x_1^2;x_2^2)$. This example, motivate us to
extend the notion of affine fractals to projective fractals. Let
$\Proj^n(\Q )$ denote the projective space of dimension $n$ and
$f_i=(\phi_1,...,\phi_n)$ for $i=1$ to $n$ denote endomorphisms
which consist of $n$ homogeneous polynomials of degree $m_i$.
Then, one can construct projetive fractals inside $\Proj^n(\bar
{\Q} )$ with respect to these homogeneous similarity maps. The
whole $\Proj^n(\Q)$ is self-similar but not a fractal, since it is
disjoint union of infinitely many copies of itself.

% ---------------------------------------------------------------
\section{Fractals in arithmetic geometry}

In general, there is no global norm on the set of points in a
fractal to motivate us how to define the notion of
fractal-dimension. In special cases, arithmetic height functions
are appropriate replacements for the norm of an algebraic
integer, particularly because finiteness theorems hold in this
context.

Northcott associated a heights function to points on the
projective space which are defined over number fields [No]. In
course of his argument for the fact that, the number of periodic
points of an endomorphism of a projective space which are defined
over a given number-field  are finite, he proved that the number
of points of bounded height is finite. Therefore, one can study
the asymptotic behavior of rational points on a fractal hosted by
a projective variety. Let us formulate a general definition.
\begin{defn}
Let $V$ be a projective variety defined over a finitely generated
field $K$ and let $f_i$ for $ i=1$ to $n$ denote finite
surjective endomorphisms of $V$ defined over $K$, which are of
degrees $>1$. A subset $F\subset V(K )$ is called a self-similar fractal
with respect to $f_i$, if $F$ is almost disjoint union of its
images under the endomorphisms $f_i$, i.e. $F=\sqcup_i f_i(F)$. $F$ is called a fractal if 
$F=\cup_i f_i(F)$
\end{defn}
$\textbf{1.}$ Let $f_i$ for $i=1,...,n$ denote homogeneous
endomorphisms of a projective space defined over a global field
$K$ with each homogeneous component of degree $m_i$. Let
$F\subseteq \Proj^n(K)$ be a fractal with respect to $f_i$:
$F=\sqcup_i f_i(F)$. One can define the fractal-dimension of $F$
to be the real number $s$ for which $ \sum_i m_i^{-s}=1 $.

\begin{prop}
In the context of projective spaces,
dimension of a self-similar fractal $F$ is well-defined and well-bahaved with respect to
fractal embeddings.
\end{prop}

\textbf{Proof:}
Indeed, for the number-field case, we use the logarithmic height
$h$ to control the height growth of points under endomorphisms.
Again we claim that if $\sum_i m_i^{-s}<1$ and $F\subseteq\cup_i
f_i(F)$ then the number of elements of $F$ of logarithmic height
less than $x$, which we denote again by $N(x)$, is bounded above
by $cx^s$ for some constant $c$ and large $x$. Let $F_i=f_i(F)$,
and $N_i(x)$ denote the number of elements of $F_i$ of logarithmic
height less than $x$. We have
$$
N(x)\leq \sum_i N_i(x)
$$
and for $f \in F_i$ and $f_i^{-1}(f)\in F$ we have
$h(f_i(f))=m_i.h(f)+O(1)$. Therefore
$$
N(x)\leq \sum_i N(m_i^{-1}x+t)
$$
for some $t$. We define a function $\bar h:[1,\infty ] \to \Real$
by $\bar h(x)=x^{-s}N(x)$. The argument of theorem 1.3 implies that
$\bar h$ is bounded, and hence the claim follows. By a similar
argument, if $F\supseteq\sqcup_i f_i(F)$ and if $r$ is a real
number such that $\sum_i m_i^{-r}>1$ then $N(x)$ is bounded below
by $cx^s$ for some constant $c$ and large $x$. One can follow the
argument of proposition 1.5 to finish the proof. $\square$

$\textbf{2.}$
For the function field case, one could use another appropriate
height function. Let $\F_q(X)$ denote the function field of an
absolutely irreducible projective variety $X$ which is
non-singular in codimension one, defined over a finite field
$\F_q$ of characteristic $p$. One can use the logarithmic height
on $\Proj^n(\F_q(X))$ defined by Neron [La-Ne]. Finiteness
theorem holds for this height function as well.

$\textbf{3.}$
Let $h,R,w,r_1 ,r_2 ,d_K ,\zeta_K$ denote class number,
regulator, number of roots of unity, number of real and complex
embeddings, absolute discriminant and the zeta function
associated to the number field $K$. Schanuel proved that [Scha]
the asymptotic behavior of points in $\Proj^n(K)$ of logarithmic
height bounded by $log(x)$ is given by
$$
{hR \over {w\zeta_K(n+1)}}\left( {2^{r_1}(2\pi)^{r_2} \over
{d_K^{1/2}}}\right)^{n+1}(n+1)^{r_1+r_2-1}x^{n+1}.
$$
This proves that rational points on projective space can not be
regarded as a fractal of finite dimension. 

$\textbf{4.}$
Schmidt in case $K=\Q$ [Schm] and Thunder for general number
field $K$ [Th] generalized the estimate of Shanuel to Grassmanian
varieties, and proved that
$$
C(G(m,n)(K),log(x))\sim c_{m,n,K}x^n
$$
where $C$ denotes the number of points of bounded logarithmic
height and $c_{m,n,K}$ is an explicitly given constant. Also,
Franke-Manin-Tschinkel provided a generalization to flag
manifolds [Fra-Man-Tsh]. Let $G$ be a semi-simple algebraic group
over $K$ and $P$ a parabolic subgroup and $V=P\backslash G$ the
associated flag manifold. Choose an embedding of
$V\subset\Proj^n$ such that the hyperplane section $H$ is
linearly equivalent to $-sK_V$ for some positive integer $s$,
then there exists an integer $t\geq 0$ and a constant $c_V$ such
that
$$
C(V(K),x)^s=c_Vx(log x)^t.
$$
All of these spaces are self-similar objects which have the
potential to be ambient spaces for fractals, but they are too huge
to be fractals themselves.  

$\textbf{5.}$ Wan proved that [Wan]
in the function field case, the asymptotic behavior of points in
$\Proj^n(K)$ of logarithmic height bounded by $d$ is given by
$$
{{hq^{(n+1)(1-g)}} \over {(q-1)\zeta_X(n+1)}}q^{(n+1)d}.
$$
which shows that $\Proj^n(\F_q(X))$ can indeed be considered as a
finite dimensional fractal. 

$\textbf{6.}$ Let $A$ be an
abelian variety over a number-field $K$ and let $F\subseteq
A(\bar {\Q})$ be a fractal with respect to endomorphisms $\phi_i$
which are translations of multiplication maps $[n_i]$ by elements
of $A(\bar {\Q})$. We define dimension of $F$ to be the real
number $s$ for which $ \sum_i n_i^{-s}=1 $. Then dimension of $F$
is well-defined and well-bahaved with respect to fractal
embeddings.
\begin{prop}
In the contect of abelian varieties,
dimension of a self-similar fractal $F$ is well-defined and well-bahaved with respect to
fractal embeddings.
\end{prop}

\textbf{Proof:}
Indeed, in this case, we use the Neron-Tate logarithmic height
$\hat h$ to control the growth of the heights of points under the
action of endomorphisms $\phi_i$. The same proof as before works
except that
$$
\hat h([n_i](f))=(n_i)^2\hat h (f)
$$
does not hold for translations of the form $[n_i]$. One should use
the fact that for the N\'{e}ron-Tate height associated to a
symmetric ample bundle on $A$ and for every $a\in A(\bar {\Q})$
and $n\in \N$, we have
$$
\hat h([n](f)+a)+\hat h([n](f)-a)=2\hat h([n](f))+ 2\hat h(a).
$$
This helps to get the right estimate. The rest of proof goes as
before.$\square$

The above notion of dimension implies that the number of points
of bounded height defined over a fixed number-field has
polynomial growth, which gives an immediate proof for 
a classical result of N\'{e}ron [Ner].

$\textbf{7.}$ Analogous to abelian varieties, one also can define
fractals on $t$-modules. By a $t$-module of dimension $N$ and
rank $d$ defined over the algebraic closure $\bar k=\overline
{\F_q(t)}$ we mean, fixing an additive group $(\G_a)^N(\bar k)$
and an injective homomorphism $\Phi$ from the ring $\F_q[t]$ to
the endomorphism ring of $(\G_a)^N$ which satisfies
$$
\Phi(t)=a_0F^0+...+a_dF^d
$$
with $a_d$ non-zero, where $a_i$ are $N\times N$ matrices with
coefficients in $\bar k$, and $F$ is a Frobenius endomorphism on
$(\G_a)^N$. One can think of polynomials $P_i\in \F_q[t]$ of
degrees $r_i$ for $ i=1$ to $n$ as self-similarities of the
$t$-module $(\G_a)^N$ and let $F\subseteq (\G_a)^N(\bar k)$ be a
fractal with respect to $P_i$ , i.e. $F=\sqcup_i \Phi(P_i)(F)$.
We define the fractal dimension of $F$ to be the real number $s$
such that $\sum_i (r_id)^{-s}=1$. Then dimension of $F$ is
well-defined and well-bahaved with respect to inclusions.
\begin{prop}
In the context of $t$-modules
dimension of a self-similar fractal $F$ is well-defined and well-bahaved with respect to
fractal embeddings.
\end{prop}

\textbf{Proof:}
Indeed, Denis defines a canonical height $\hat h$ on $t$-modules
which satisfies
$$
\hat h[\Phi(P)(\alpha)]=q^{dr}.\hat h[\alpha]
$$
for all $\alpha \in (\G_a)^N$, where $P$ is a polynomial in
$\F_q[t]$ of degree $r$ [Den]. This can be used to prove the result
in the same lines as before. One can get information on the
asymptotic behavior of $N(\G_a^N(\bar k),x)$ by representing
$\G_a^N(\bar k)$ as a fractal.$\square$

% ---------------------------------------------------------------
\section{Diophantine approximation by fractals}

This section is devoted to proving theorems which were mentioned
in the introduction. The arguments are along the same lines as
analogous classical results.

Roth's theorem on Diophantine approximation of rational points on
projective line implies a version on projective varieties defined
over number-fields. Self-similarity of rational points on abelian
varieties makes room to improve the estimates. This argument can
be imitated in case of arithmetic fractals defined over finitely generated fileds.

\begin{thm} (Fractal version of Roth's thereom on diopphantine approximation)
Fix a finitely generated field of characteristic zero $K$ 
and $\sigma :K\hookrightarrow \Cplx$ a
complex embedding. Let $V$ be a smooth projective algebraic
variety defined over $K$ and let $L$ be an very ample line-bundle on
$V$. Denote the arithmetic height function associated to the
line-bundle $L$ by $h_L$. Suppose $F\subset V(K)$ is a fractal
subset with respect to finitely many height-increasing
self-endomorphisms $\phi_i:V\to V$ defined over
$K$ such that for
all $i$ we have
$$
h_L(\phi_i(f))>m_ih_L(f)+0(1)+
$$
where $m_i>1$. Fix a Riemannian metric on $V_{\sigma}(\Cplx)$ and
let $d_{\sigma}$ denote the induced metric on
$V_{\sigma}(\Cplx)$. 
Then, for every $\delta>0$ and every choice
of an algebraic point $\alpha\in V(\bar {K})$ which is not a
critical value of any of the $\phi_i$'s and all choices of a
constant $C$, there are only finitely many fractal points
$\omega\in F$ approximating $\alpha$ such that 
$$
d_{\sigma}(\alpha ,\omega)\leq Ce^{-\delta h_L(\omega)}.
$$
\end{thm}

\begin{prop}
With assumptions of the above theorem, suppose for some $\delta_0>0$ 
we have that, for any choice of a constant $C$
and every choice
of an algebraic point $\alpha\in V(\bar {K})$
there are only finitely many fractal points
$\omega\in F$ approximating $\alpha$ in the following manner
$$
d_{\sigma}(\alpha ,\omega)\leq Ce^{-\delta_0 h_L(\omega)}.
$$
Then, for every $\delta>0$ and every choice
of an algebraic point $\alpha\in V(\bar {K})$ which is not a
critical value of any of the $\phi_i$'s and all choices of a
constant $C$, there are only finitely many fractal points
$\omega\in F$ approximating $\alpha$ such that 
$$
d_{\sigma}(\alpha ,\omega)\leq Ce^{-\delta h_L(\omega)}.
$$
\end{prop}

\textbf{Proof (Proposition).} Note that, we have assumed that the above is true for some
$\delta_0>0$ without any assumption on $\phi_i$ or on $\alpha$.
Let $\delta'>0$ be infimum of such $\delta_0>0$.

Fix $\epsilon>0$ such that $\epsilon<\delta' <m_i\epsilon$ for
all $i$. Suppose that $w_n$ is an infinite sequence of elements in
$F$ such that $\omega_n\to \alpha$ which satisfies the estimate
$$
d_{\sigma}(\alpha ,\omega_n)\leq Ce^{-\epsilon h_L(\omega_n)}.
$$
then infinitely many of them are images of elements of $F$ under
the same $\phi_i$. By going to a subsequence, one can find a
sequence $\omega'_n$ in $F$ and an algebraic point $\alpha'$ in
$V(\bar {K})$ such that $\omega'_n \to \alpha'$ and for a fixed
$\phi_i$ we have $\phi_i(\alpha')=\alpha$ and
$\phi_i(\omega'_n)=\omega_n$ for all $n$. Then
$$
d_{\sigma}(\alpha ,\omega_n)\leq Ce^{-\epsilon h_L(\omega_n)}\leq
C'e^{-\epsilon m_i h_L(\omega'_n)}
$$
for an appropriate constant $C'$. On the other hand,
$$
d_{\sigma}(\alpha' ,\omega'_n)\leq C''d_{\sigma}(\alpha ,\omega_n)
$$
holds for an appropriate constant $C''$ and large $n$ by
injectivity of $d\phi_i^{-1}$ on the tangent space of $\alpha$.
This contradicts our assumption on $\delta'$, because $\delta' <m_i\epsilon$.
$\square$

\textbf{Proof (Theorem).}  If we assume that points of $F$ and similarity maps are defined over some
number-field, Roth's theorem implies that the assumption of theorem is true for any 
$\delta_0>2$. All such examples are forward orbits of finitely many 
height increasing self-similarities. The same is true for finitely generated field of characteristic zero
by a result of Lang [Lan] generalizing Roth's theorem and height defined by Moriwaki [Mor].$\square$

\begin{rmk}
The conditions of thereom could not hold true for general fractals in
$V(\bar {K})$. For example, torsion points of an abelian variety
are dense in complex topology, and have vanishing height.
Therefore, our fractal analogue of Roth's theorem could not hold
in this case.
\end{rmk}

Let us state a more precise version of our version of Siegel's theorem.

\begin{thm}(Fractal version of Siegel's theorem on integral points) Fix a
finitely generated field of characteristic zero $K$. Let $V$ be a smooth affine algebraic variety
defined over $K$ with smooth projectivization $\bar V$ and let $L$
be an very ample line-bundle on $\bar V$. Denote the arithmetic height
function associated to the line-bundle $L$ by $h_L$. Suppose
$F\subset V({K})$ is a fractal subset with respect to
finitely many height-increasing polynomial self-endomorphisms
$\phi_i:V\to V$ defined over $K$ such that for all $i$ we have
$$
h_L(\phi_i(f)) > m_ih_L(f)+0(1)
$$
where $m_i>1$. One could also replace this assumption with norm
analogue. For any affine hyperbolic algebraic curve $X$ embedded
in $V$ defined over $K$ we have $X(K)\cap F$ is a
finite set.
\end{thm}

We borrow a lemma from [Ser] whose proof goes exactly as in that reference.

\begin{lem}
Let $K$ be a finitely generated field of characteristic zero. Let $X$ be a curve defined over $K$.
Assume genus of $X$ is $\geq 1$.
If $P_n$ is a sequence of distince
points in $X(K)$, which means that their heights tends to infinity and if
$\phi$ defined over $K$ is a non-constant rational function on $X$. From some point
on, no $P_n$ is pole of $\phi$. 
Then for $z_n=\phi(P_n)$ which a point of the
projective space defined over $K$ we have
$$
\lim_{n\to \infty} {{log|z_n|_{v}}\over {log H(z_n)}}=0
$$
\end{lem}

\textbf {Proof.}
Assume this is false. By taking a subsequence and replacing $\phi$ by
$1/\phi$, we may suppose that
$$
{{log|z_n|_{v}} \over {log H(z_n)}}\rightarrow \lambda
$$
where $-\infty <\lambda< 0$.
In particular, $z_n \rightarrow 0$ in $K_v$ and by taking a subsequence, we may assume
that $P_n$ converges to a zero $P_0$ of $\phi$. As we are on a curve, $P_0$ is an algebraic
point of X.
Between $H(P_n )$, the height corresponding to a morphism $X \rightarrow \mathbb P_N$,
and $H(z_n )$, corresponding to a morphism $X \rightarrow \mathbb P_1$, we have an inequality
$$
H(z_n )\ll H(P_n )^{l}
$$
for some positive $l$.
On the other hand, if $e$ is the multiplicity of $P_0$ as a zero of $\phi$, we have
$|z_n|_{v}\approx d_{v}(P_n, P_0)^{e}$.
Therefore, there is $c > 0$ such that for sufficiently large $n$,
$$
d_{v}(P_n, P_0)\leq 1/H(P_n )^c
$$
which contradicts the approximation theorem.$\square$

\textbf {Sketch of Proof of Theorem.} Let $\sigma :K\hookrightarrow \Cplx$
denote a complex embedding of $K$. Fix a Riemannian metric on
$V_{\sigma}(\Cplx)$ and let $d_{\sigma}$ denote the induced
metric on $V_{\sigma}(\Cplx)$. Then by our version of Roth's
theorem, for every $\delta>0$ and every choice of an algebraic
point $\alpha\in V(\bar {K})$ which is not a critical value of
any of the $\phi_i$'s and all choices of a constant $C$, there
are only finitely many fractal points $\omega\in F$ approximating
$\alpha$ such that
$$
d_{\sigma}(\alpha ,\omega)\leq H_L(\omega)^{-\delta }.
$$
where $log(H_L)=h_L$. 

In case $K$ is trancendental, we have to pick a model for $K$ over algebraic closure of 
$\mathbb{Q}$ in $K$ following Lang [lan].
Now if $P_n$ is a sequence of distince
points in $X(K)\cap F$, their heights tends to infinity and if
$\phi$ is a non-constant rational function on $X$ from some point
on no $P_n$ is pole of $\phi$.
Then by above proposition
$$
\lim_{n\to \infty} {{log|z_n|_{\sigma}}\over {log H(z_n)}}=0
$$
On the other hand, one defines height of rational points by
$$
H(z)=\prod_{v\in M_K} sup(1,|z|_v),
$$
where $|.|_v$ are normalized according to a product formula.
Since similarity maps of $F$ are expanding, we know 
that $F$ is forward orbit of finitely many points. So for a
finite set of places $S$ we have
$$
H(z)=\prod_{v\in S} sup(1,|z|_v),
$$
and therefore
$$
log H(z)=\sum_{v\in S} log(sup(1,|z|_v)).
$$
Then, we have 
$$
1=\sum_{v\in S} sup(0,{{log|z_n|_{\sigma}}\over {log H(z_n)}})\leq \sum_{v\in S}  {{log|z_n|_{\sigma}}\over {log H(z_n)}}
$$
which could not be true, because the above limit is zero. This
implies the finiteness result we are seeking for. $\square$

\begin{rmk} If $F$ and its self-similarity maps are defined over $K$, 
then $F$ is forward orbit of finitely many points which are not neccessarily algebraic, 
and the above result 
is not implied by Siegel's theorem for $S$-integral points.
\end{rmk}

This being true, we expect the following version of Liouville's theorem on diophantine approximation holds:

\begin{thm}(Fractal version of Liouville's theorem on diophantine approximation)
Fix a finitely generated field of characteristic zero $K$. Let $V$ be a smooth projective algebraic variety
defined over $K$ and let $L$
be an very ample line-bundle on $V$. Denote the arithmetic height
function associated to the line-bundle $L$ by $h_L$. 
Then there exists a positive constant $\delta_0$ such that for any positive constant $c$
For any geometrically irreducible algebraic subvariety $E$ of $V$ defined over $K$ 
and $d_w(x,E)$ denoting the $w$-adic distance from $x$ to $E$,
there are only finitely many points defined over $K$ in $V(K)$ outside $E$ satisfying the following inequality
$$
d_w(x,E)< cH(x)^{-\delta_0}
$$
except for points in an algebraic variety $V(\delta_0 )$ which is of strictly smaller dimension of $V$.
\end{thm}

\textbf {Proof.} This is a weak form of Vojta conjectures.
In the number field case, this is mentioned in Faltings-Wustholz [Fa-Wu] as a trivial result in case $E$ is geometrically irreducible.
In the case of finitely generated fields of characteristic zero the result is a consequence of theorem I' in seminal work of Lang [Lan].$\square$

Now, the following version of Falting's theorem, can be proved 
using the methods of self-similarity and height expansion.

\begin{thm}(Fractal version of Faltings' theorem on diophantine approximation on abelian varieties)
Fix a finitely generated field of characteristic zero $K$. Let $V$ be a smooth projective algebraic variety
defined over $K$ and let $L$
be an very ample line-bundle on $V$. Denote the arithmetic height
function associated to the line-bundle $L$ by $h_L$. Suppose
$F\subset V(K)$ is a fractal subset with respect to
finitely many height-increasing polynomial finite self-endomorphisms
$\phi_i:V\to V$ defined over $K$ such that for all $i$ we have
$$
h_L(\phi_i(f)) > m_ih_L(f)+0(1)
$$
where $m_i>1$. Fix any positive constants $\kappa$ and $c$. 
For any irreducible algebraic subvariety $E$ of $V$ defined over $K$ 
and $d_w(x,E)$ denoting the $w$-adic distance from $x$ to $E$,
there are only finitely many points defined over $K$ in $F$ outside $E$ satisfying the following inequality
$$
d_w(x,E)< cH(x)^{-\kappa}
$$
\end{thm}

\textbf {Sketch of Proof.}
Reduction of the inequalty for some positive $\delta_0$ to arbitrary $\delta >0$ is done in the 
same manner as in our version of Roth's theorem. Getting rid of $V(\delta )$ is the result of the fact that
$V(\delta )$ is invariant under $\phi_i$ and $F\cap V(\delta )$ is again a fractal. One can proceed by reducing the problem from $V$ and $E$
to $V(\delta )$ and $E\cap V(\delta )$ and applying induction.$\square$

The following will be a special case:

\begin{cor}
Fix a number field $K$. Let $E$ be an irreducible affine smooth algebraic variety
defined over $K$. If
$d_w(x,E)$ denotes the $w$-adic distance from $x$ to $E$, 
then for any positive constant $\delta$ such that for any positive constant $c$
there are only finitely many points defined over ring of integers $O_K$ outside $E$ satisfying the following inequality
$$
d_w(x,E)< cH(x)^{-\delta}
$$
\end{cor}
In particular, we have the following be true:

\begin{cor}
Let $D$ be an irreducible affine smooth divisor
defined over $\mathbb Q$. If
$d_w(x,E)$ denotes the $w$-adic distance from $x$ to $E$, 
then for any positive constant $\delta$ such that for any positive constant $c$
there are only finitely many points defined defined over $\mathbb Z$ outside $D$ satisfying the following inequality
$$
d(x,E)< c||x||^{-\delta}
$$
\end{cor}

A simple highschool implication would be the following 

\begin{cor}
Let $f$ be an algebraic equation in two variables determining an irreducible algebraic curve $C$ in $\mathbb R^2$.
Then for positive constant $\delta$ and for any positive constant $c$
there are only finitely many points defined in $\mathbb Z^2$ outside the curve $C$ satisfying the following inequality
$$
d(x,C)< c||x||^{-\delta}
$$
\end{cor}

% ---------------------------------------------------------------
\section{Fractal conjecture}
\begin{conj} (Fractal conjecture)
Let $V$ be an irreducible variety defined over a finitely
generated field $K$ and let $F\subset V(K)$ denote a fractal
on $V$ with respect to finitely many height-increasing self-maps
$$
f_i:V \rightarrow V 
$$ 
defined over $K$.
Then, for any reduced subscheme $Z$ of $V$ defined over
$K$ the Zariski closure of $Z(\bar K)\cap F$ is union of finitely
many points and finitely many components $B_j$ such that $B_j(K)\cap F$ 
is a fractal in $B_j$ for each $j$, with respect to some of $f_i$.
If non of the components of $B_j$ are pre-priodic with respect to any of $f_i$
then any $B_j(K)\cap F$ 
is a fractal in $B_j$ with respect to all of $f_i$.
\end{conj}

\begin{rmk}
You can start with $F\subset V(\bar K)$, but then you can not assume 
$f_i$ are height increasing and instead you may join some $B_j$ to make a fractal.
\end{rmk}

In particular, we have stated the following 
 
\begin{conj}
For any algebraic curve $C$ embedded
in $V$ defined over $K$ which is not invarient under $f_i$, 
we have $C(\bar K)\cap F$ is at most a finite set.
\end{conj}

The following would be a corollary

\begin{cor} 
For any hyperbolic projective curve $C$ embedded in an abelian variety $A$
and any finitely generated subgroup $\Gamma$ in $A$ all af them 
defined over a finitely generated field of characteristic zero $K$, we have 
$C(\bar K)\cap \Gamma$ is finite. Even $C(\bar K)\cap Div(\Gamma)$ is finite, where $Div(\Gamma)$ is divisible group of $\Gamma$.
\end{cor}

\begin{rmk}
In case $A$, $X$ and $\Gamma$ are defined on a number field, the above is content of
a conjecture of Mordell, proved by G. Faltings [Fal1].
\end{rmk}

There is also another implication of our self-similarity conjecture:

\begin{conj}
(Forward orbit conjecture) Let $V$ be an irreducible variety defined over a finitely
generated field $K$ and let $f_i:V \rightarrow V$ denote finitely many self maps of $V$ defined over $K$. 
Let $F$ denote the forward orbit with respect to $f_i$ of finitely many points of $V$ defined over $K$. 
Then, for any reduced subscheme $Z$ of $V$ defined over
$K$ the Zariski closure of $Z(\bar K)\cap F$ is union of finitely
many points and finitely many components $B_j$ such that $B_j(\bar K)\cap F$ 
is the forward orbit with respect to some $f_i$ of finitely many points of $B_j$ defined over $K$, for each $j$.
\end{conj}

It is instructive to notice that, the common geometric structures
appearing in the context of Diophantine geometry, is exactly the
same as the objects appearing in dynamics of endomorphisms of
algebraic varieties which was the original context that height
functions were introduced.

Let us start by restating Raynaud's theorem on torsion points of
abelian varieties lying on a subvariety [Ray], which is a special case of
conjecture 5.1.
\begin{thm}
(Raynaud) Let $K$ be a number field and let $A$ be an abelian variety over the algebraically
closed field $\bar K$, and $Z$ a reduced
subscheme of $A$. Then every irreducible component of the Zariski
closure of $Z(\bar K)\cap A(\bar K)_{tor}$ is a translation of an
abelian subvariety of $A$ by a torsion point.
\end{thm}

Another special case is Faltings' theorem on finitely generated
subgroups of abelian varieties which has a very similar feature
[Fal].
\begin{thm}
(Faltings) Let $K$ be a number field and let $A$ be an abelian variety over the algebraically
closed field $\bar K$,
 and $\Gamma$ be a
finitely generated subgroup of $A(\bar K)$. For a reduced
subscheme $Z$ of $A$, every irreducible component of the Zariski
closure of $Z(\bar K)\cap \Gamma$ is a translation of an abelian
subvariety of $A$.
\end{thm}

Another consequence of conjecture 05.1 would be the following
version of generalized Lang's conjecture [Zha].

\begin{conj}
(S. Zhang) Let $X$ be an algebraic variety defined over a number-field $K$
and let $f:X \to X$ be a surjective endomorphism defined over
$K$. Suppose that the subvariety $Y$ of $X$ is not pre-periodic in
the sense that the orbit $\{Y,f(Y),f^2(Y),...\}$ is not finite,
then the set of pre-periodic points in $Y$ is not Zariski-dense
in $Y$.
\end{conj}

Lang's conjecture is confirmed by Raynaud's result mentioned
above in the case of abelian varieties and by results of Laurent
[Lau] and Sarnak [Sar-Ada-Rud] and S. Zhang [Zha] in the case of multiplicative
groups.

%%%%%%%%%%%%%%%%%%%%%%%%%%%%%%%%%%%%%%%%%%%%%%%%%%%%%%%%%%%%%%%%%
\section{Quasi-fractal conjecture}

Andre-Oort conjecture on sub-varieties of Shimura varieties is
motivated by conjectures of Lang and Manin-mumford which were
proved by Raynaud and Faltings as mentioned above. Motivated by
the Andre-Oort conjecture (look at [Ed] for an exposition of this
conjecture), we also present another conjecture in the same lines
for quasi-fractals in an algebraic variety $X$, where
self-similarities are allowed to be induced by geometric
self-correspondences on $X$ instead of self-maps. For
quasi-fractals, we drop the requirement that similar images shall
be almost-disjoint.
\begin{conj}
(Quasi-fractal conjecture)
Let $V$ be an irreducible variety defined over a finitely
generated field $K$ and let $F\subset V(\bar K)$ denote a
quasi-fractal on $V$ with respect to correspondences $Y_1,...,Y_n$
on $V$ living in $V\times V$ with both projections finite and
surjective. $F$ may contain
a subvariety of $V$. Then, for any reduced subscheme $Z$ of $V$ defined
over $K$ the Zariski closure of $Z(\bar K)\cap F$ is union of
finitely many points and finitely many components $B_i$ such that
for each $i$ the intersection $B_i(\bar K)\cap F$ is a
quasi-fractal in $B_i$ with respect to some correspondences induced by
$Y_i$.
\end{conj}

Conjecture 6.1 is a more sophisticated version of our previous
conjecture, which implies fractal conjecture also absorbs Andre-Oort conjecture into the
fractal formalism. This version utilizes the concept of
quasi-fractals.
\begin{defn}
Let $V$ be an algebraic variety and let $Y_i\hookrightarrow
V\times V$ for $ i=1$ to $n$ denote correspondences on $V$ where
$\pi_1$ and $\pi_2$ are projection to the first and second factor
in $V\times V$ which are finite and surjective when restricted on
the image of each $Y_i$ in $V\times V$ for $ i=1$ to $n$. A
subset $F\subseteq V(\bar K)$ is called a quasi-fractal with respect to
$Y_1,...,Y_n$ if $F$ is invariant under the action of
correspondences $Y_1,...,Y_n$.
\end{defn}

The $l$-Hecke orbit of a point on the moduli-space of principally
polarized abelian varieties is an example of a quasi-fractal with
respect to the $l$-Hecke correspondences associated to
$l$-isogenies.

J. Pila also gives an unconditional proof of the Andre-Oort conjecture for arbitrary 
products of modular curves [Pil]. He with cooperation Tsimerman eventually proves the full Andre-Oort conjecture [Pil-Tsi1] [Pil-Tsi2].

Conjecture 6.1 also covers a parallel version of Andre-Oort
conjecture for $l$-Hecke orbit of a special point in the function
field case [Bre].

%---------------------------------------------------------------
\subsection*{acknowledgements}
I have benefited from conversations with N. Fakhruddin,M. Hadian, H.
Mahdavifar, A. Nair, O. Naghshineh, A. Rajaei, P. Sarnak, C. Soule, V. Srinivas , N. Talebizadeh, P. Vojta
for which I am thankful. Peter Sarnak particularly gave crucial comments which led 
to the final version of the paper. 
I would also like to thank Sharif University
of Technology and Young Scholars Club who partially supported this
research and the warm hospitality of TIFR, ICTP, IHES and Princeton University where
this work was written, completed and revised.
%-----------------------------------------------------------------

%-------------------------------------------------------------------------

Sharif University of Technology, e-mail: rastegar@sharif.edu
\\Princeton University, e-mail: rastegar@princeton.edu

\end{document}